\documentclass[11pt]{article}
\usepackage{latexsym,amsmath,amssymb,amscd}
\begin{document}

\makeatletter
\@addtoreset{figure}{section}
\def\thefigure{\thesection.\@arabic\c@figure}
\def\fps@figure{h,t}
\@addtoreset{table}{bsection}

\def\thetable{\thesection.\@arabic\c@table}
\def\fps@table{h, t}
\@addtoreset{equation}{section}
\def\theequation{%\thesection.
\arabic{equation}}
\makeatother

\newcommand{\bfi}{\bfseries\itshape}

\newtheorem{theorem}{Theorem}
\newtheorem{acknowledgment}[theorem]{Acknowledgment}
\newtheorem{algorithm}[theorem]{Algorithm}
\newtheorem{axiom}[theorem]{Axiom}
\newtheorem{case}[theorem]{Case}
\newtheorem{claim}[theorem]{Claim}
\newtheorem{conclusion}[theorem]{Conclusion}
\newtheorem{condition}[theorem]{Condition}
\newtheorem{conjecture}[theorem]{Conjecture}
\newtheorem{construction}[theorem]{Construction}
\newtheorem{corollary}[theorem]{Corollary}
\newtheorem{criterion}[theorem]{Criterion}
\newtheorem{definition}[theorem]{Definition}
\newtheorem{example}[theorem]{Example}
\newtheorem{lemma}[theorem]{Lemma}
\newtheorem{notation}[theorem]{Notation}
\newtheorem{problem}[theorem]{Problem}
\newtheorem{proposition}[theorem]{Proposition}
\newtheorem{remark}[theorem]{Remark}
\numberwithin{theorem}{section}
\numberwithin{equation}{section}

%%% Todo
\newcommand{\todo}[1]{\vspace{5 mm}\par \noindent
\framebox{\begin{minipage}[c]{0.95 \textwidth}
\tt #1 \end{minipage}}\vspace{5 mm}\par}
%%%

\newcommand{\1}{{\bf 1}}

\newcommand{\ev}{{\rm ev}}
\newcommand{\id}{{\rm id}}
\newcommand{\ie}{{\rm i}}
\newcommand{\Hom}{{\rm Hom}\,}
\newcommand{\Ker}{{\rm Ker}\,}
\newcommand{\lf}{{\rm l}}
\newcommand{\Ran}{{\rm Ran}\,}
\newcommand{\Tr}{{\rm Tr}\,}

\newcommand{\G}{{\rm G}}
\newcommand{\U}{{\rm U}}

\newcommand{\Bc}{{\mathcal B}}
\newcommand{\Cc}{{\mathcal C}}
\newcommand{\Hc}{{\mathcal H}}
\newcommand{\Oc}{{\mathcal O}}

\pagestyle{myheadings}
\markboth{\sl Belti\c t\u a and Ratiu: Geometric
representations}
{\sl Belti\c t\u a and Ratiu: Geometric
representations}

%\begin{document}

\makeatletter
\title{\textbf{Geometric representation theory for
unitary groups of operator algebras}}
\author{Daniel Belti\c t\u a$^{1}$ and  Tudor S.
Ratiu$^{2}$}
\addtocounter{footnote}{1}
\footnotetext{Institute of Mathematics ``Simion
Stoilow'' of the Romanian Academy, 
014700 Bucharest, Romania.
\texttt{Daniel.Beltita@imar.ro}. }
\addtocounter{footnote}{1}
\footnotetext{Section de Math\'ematiques and Centre  Bernoulli, 
\'Ecole Polytechnique F\'ed\'erale de Lausanne, 
CH-1015 Lausanne,
Switzerland. \texttt{Tudor.Ratiu@epfl.ch}.}
%\date{November 26, 2004}
\makeatother
\maketitle

\begin{abstract}
Geometric realizations for the
restrictions of GNS representations to unitary groups
of $C^*$-algebras are constructed. 
These geometric realizations use an
appropriate concept of reproducing kernels on vector
bundles. To build such realizations in spaces of
holomorphic sections, a class of
complex coadjoint orbits of 
the corresponding real Banach-Lie groups are described  
and some homogeneous holomorphic
Hermitian vector bundles that are naturally associated
with the coadjoint orbits are constructed. 


\end{abstract}

\section{Introduction}

The study of geometric properties of state spaces is a
basic topic in  the theory of operator algebras (see,
e.g.,  \cite{AS01} and \cite{AS03}).
The GNS construction produces representations of
operator algebras out of states. 
From this point of view, we think it interesting to
investigate the geometry behind these
representations. 

One method to do this is to proceed as in the theory 
of geometric realizations of Lie group representations 
(see e.g., \cite{DF95}, \cite{Do97}, \cite{Ne00},
\cite{ES02}) and to try to build the representation
spaces as spaces of sections of certain vector bundles. 
The basic ingredient in this construction is the
reproducing kernel Hilbert space 
(see, e.g., \cite{Ar50}, \cite{Sch64},
\cite{Ha82}, \cite{Sa88}, \cite{Arv98},
\cite{BH98}, \cite{Ne00}).

In the present paper we show that the aforementioned
method  can indeed be applied to the case of group
representations obtained by restricting GNS
representations to unitary groups of $C^*$-algebras. 
More precisely, for these
representations, 
we construct one-to-one intertwining operators from
the representation spaces onto reproducing kernel
Hilbert spaces of sections of certain Hermitian vector
bundles (Theorem~\ref{Theorem12}). 
The construction of these vector bundles is based on a
choice of a 
sub-$C^*$-algebra that is related in a suitable way to
the state involved in the GNS construction (see
Construction~\ref{homog}). 
It turns out that, in the case of normal states of
$W^*$-algebras, there is a natural choice of the
subalgebra 
(namely the centralizer subalgebra), 
and the base of the corresponding vector bundle is
just one of the symplectic leaves 
studied in our previous paper~\cite{BR04}.  
Since the corresponding symplectic leaves are just
unitary orbits of states, the geometric representation
theory initiated in the present paper provides, in 
particular, a geometric interpretation of the
result in 
\cite{GK60}, namely the equivalence class of an
irreducible GNS representation only depends on the
unitary orbit of the corresponding pure state. 

In~\cite{Boy93} and references therein one can find 
several interesting results regarding the
classification of unitary group representations of
various operator algebras.  The point of the present
paper is to show that some of these representations 
(namely the ones obtained by restricting GNS
representations to unitary groups) can be realized
geometrically following the pattern of the classical
Borel-Weil theorem for compact groups. This raises
the challenging problem of finding
geometric realizations of more general representations
of unitary groups of operator algebras. 
Similar results for other classes of
infinite-dimensional groups have been already
obtained: 
see \cite{NRW01}, \cite{DPW02},~\cite{Wo04} for
direct limit groups, 
and \cite{Boy80}, \cite{Ne02},~\cite{Ne04} 
for groups related to operator ideals. 
The same problem of geometric realizations for
representations of the restricted unitary group was
raised at the end of~\cite{Boy88}.

The structure of the paper is as follows. 
Since the reproducing kernels we need in the present
paper 
show up most naturally in a $C^*$-algebraic setting 
(Construction~\ref{kernel}), we establish in
Section~2 
the appropriate versions of a number of results in
\cite{BR04}. 
Section~3 gives a general construction of
homogeneous Hermitian vector bundles associated with
GNS representations. 
Section~4 is devoted to the concept of reproducing
kernel suitable for the applications we have in mind. 
In Section~5 we construct such reproducing kernels out
of GNS representations and we prove our main theorems
on geometric realizations of GNS representations 
(Theorems \ref{Theorem12}~and~\ref{holomorphy}). 

\section{Coadjoint orbits and $C^*$-algebras with
finite traces}

In this section we extend to a $C^*$-algebraic
framework 
a number of results that were proved in \cite{BR04}
for symplectic leaves in preduals of $W^*$-algebras. 
We begin by establishing some notation that will be
used throughout the paper. 

\begin{notation}\label{Notation1}
\normalfont
For a unital $C^*$-algebra $A$ with the unit $\1$ we
shall use the following notation: 
$$\begin{aligned}
\{a\}'&=\{b\in A\mid ab=ba\}\quad \mbox{whenever }a\in
A,\\ 
A^\varphi&=\{a\in A\mid(\forall b\in
A)\quad\varphi(ab)=\varphi(ba)\}
 \quad\mbox{whenever }\varphi\in A^*,\\ 
\G_A&=\{g\in A\mid g\mbox{ invertible}\}, \\
\U_A&=\{u\in A\mid uu^*=u^*u=\1\}\subseteq\G_A,\\
A^{{\rm sa}}&=\{a\in A\mid a=a^*\},\\
(A^*)^{{\rm sa}}&=\{\varphi\in A^*\mid 
     (\forall a\in A)\quad
\varphi(a^*)=\overline{\varphi(a)}\}.
\end{aligned}$$
In the special case of a $W^*$-algebra $M$ we will
also use 
the notation 
$$\begin{aligned}
M_*&=\{\varphi\in M^*\mid \varphi\mbox{ is
$w^*$-continuous}\}, \\
M_*^{{\rm sa}}&=M_*\cap(M^*)^{{\rm sa}}. \\
\end{aligned}
$$
\end{notation}

\begin{proposition}\label{Proposition2}
Let $A$ be a unital $C^*$-algebra having a faithful tracial
state 
$\tau\colon A\to{\mathbb C}$. 
Consider the mapping
$$\Theta^\tau\colon A\to A^*,\quad
a\mapsto\Theta^\tau_a,$$
where for each $a\in A$ we define 
$$\Theta^\tau_a\colon A\to{\mathbb C},\quad
\Theta^\tau_a(b):=\tau(ab).$$
The mapping $\Theta^\tau$ has the following
properties:
\begin{itemize}
\item[{\rm(a)}] $\Ker\Theta^\tau=\{0\}$. 

\item[{\rm(b)}] 
$A^{\Theta^\tau_a}=\{a\}'$ for all $a\in A$. 

\item[{\rm(c)}] The mapping $\Theta^\tau$ is
$\G_A$-equivariant 
with respect to the adjoint action of $\G_A$ on $A$
and the coadjoint 
action of $\G_A$ on $A^*$. 
In particular, the mapping 
$$\Theta^\tau|_{A^{{\rm sa}}}\colon A^{{\rm
sa}}\to(A^*)^{{\rm sa}}$$
is $\U_A$-equivariant with respect to the adjoint
action of 
$\U_A$ on $A^{{\rm sa}}$ and the coadjoint 
action of $\U_A$ on $(A^*)^{{\rm sa}}$. 

\item[{\rm(d)}] For each $a\in A$ the mapping
$\Theta^\tau$ induces a bijection 
of the adjoint orbit $\G_A \cdot a$ onto the coadjoint
orbit $\G_A\cdot \Theta^\tau_a$. 
In particular, if $a\in A^{{\rm sa}}$ and there exists
a conditional expectation of $A$ onto $\{a\}'$ 
then we have 
a commutative diagram of $\U_A$-equivariant
diffeomorphisms of Banach manifolds
$$\begin{CD} 
\U_A \cdot a \\
@AA{\qquad\quad\searrow\hskip-1pt\Theta^\tau}A \\
U_A/\U_{\{a\}'} @>>> \U_A \cdot \Theta^\tau_a 
\end{CD}
$$
\item[{\rm(e)}] If, moreover, 
$A$ is a $W^*$-algebra and the faithful tracial state
$\tau$ is normal, then $\Ran\Theta^\tau\subseteq A_*$
and the hypothesis on conditional expectation from
{\rm(d)} holds for each $a\in A^{{\rm sa}}$. 
\end{itemize}
\end{proposition}

\noindent\textbf{Proof.\  \  } This is proved like 
Proposition~2.12 in \cite{BR04}. 
\quad $\blacksquare$

\begin{corollary}\label{Corollary3'}
Let $A$ be a unital $C^*$-algebra having a faithful tracial
state 
$\tau\colon A\to{\mathbb C}$ and let $a=a^*\in A$
be such that 
there exists a conditional expectation of $A$ onto
$\{a\}'$. 
Then the unitary orbit of $a$ has  
a natural structure of 
$\U_A$-homogeneous weakly symplectic manifold. 
\end{corollary}

\noindent\textbf{Proof.\  \  }
Use Proposition~\ref{Proposition2} along
with the reasoning that leads to Corollary~2.9 
in~\cite{OR03}. 
\quad $\blacksquare$

\begin{corollary}\label{Corollary3}
In a $W^*$-algebra $M$ that admits a
faithful normal tracial state the unitary orbit of each
self-adjoint element has a natural structure of 
$\U_M$-homogeneous weakly symplectic manifold. 
\end{corollary}

\noindent\textbf{Proof.\  \  }
Use Proposition~\ref{Proposition2} along 
with Corollary~2.9 
in ~\cite{BR04}. 
\quad $\blacksquare$

\begin{proposition}\label{Proposition4}
Let $A$ be a unital $C^*$-algebra and $a=a^*\in A$ have
the spectrum a finite set. 
Then the unitary orbit of $a$ has a natural structure
of $\U_A$-homogeneous 
complex Banach manifold. 
\end{proposition}

\noindent\textbf{Proof.\  \  }
Let $a=a^*\in A$ such that there exist finitely many
different numbers $\lambda_1,\dots\lambda_m\in{\mathbb R}$ 
with 
$$a=\lambda_1 e_1+\cdots+\lambda_m e_m,$$
where $e_1,\dots,e_m\in A$ are orthogonal projections
satisfying 
$e_ie_j=0$ whenever $i\ne j$ and $e_1+\cdots+e_m=\1$. 
It is clear that, denoting $p_j=e_1+\cdots+e_j$ for
$j=1,\dots,m-1$, we have 
$$\{a\}'%=\{b\in M\mid ab=ba\}
=\{b\in A\mid be_j=e_jb\mbox{ for }1\le j\le m\}
=\{b\in A\mid bp_j=p_jb\mbox{ for }1\le j\le m-1\}.
$$
Thus both 
the unitary orbit of $a$ 
and 
the unitary orbit of $(p_1,\dots,p_{m-1})\in
A\times\cdots\times A$ 
can be identified with 
$\U_A/\U_{\{a\}'}$, 
and now the desired conclusion follows by 
Corollary~16 in \cite{Be04} 
(see Proposition~\ref{Proposition16} below for the
special case of 
$W^*$-algebras). 
%or, alternatively, by Corollary~16 in~\cite{Be04}. 
\quad $\blacksquare$

\begin{theorem}\label{Theorem5}
If  a unital $C^*$-algebra $A$ possesses a faithful
tracial state, then
the unitary orbit of each self-adjoint element with
finite spectrum 
has a natural structure of an 
$\U_A$-homogeneous weakly K\"ahler manifold. 
\end{theorem}

\noindent\textbf{Proof.\  \  }
The proof is similar to the one of Proposition~4.8 in
\cite{BR04}.
\quad $\blacksquare$

\medskip

We conclude this section by a result that will be
needed in 
the proof of Theorem~\ref{holomorphy}. 
We say that a unital $W^*$-algebra $M$ 
is {\it finite} if for every $u\in M$ with $u^*u=\1$ we have $uu^*=\1$ as well. 
In this case each projection $p\in M$ is finite, in the sense that, 
if $v\in M$, $v^*v=p$ and $vv^*$ is a projection smaller than $p$, then $vv^*=p$. 
Clearly a unital $W^*$-algebra is finite if it has a faithful tracial state.  

\begin{proposition}\label{Proposition16}
Let $M$ be a finite $W^*$-algebra and
$e_1,\dots,e_n\in M$ orthogonal 
projections satisfying $e_ie_j=0$ whenever $i\ne j$
and 
$e_1+\cdots+e_n=\1$. 
Next consider the sub-$W^*$-algebra 
$M_0=\{a\in M\mid ae_j=e_ja\mbox{ for }j=1,\dots,m\}$ 
of $M$,
 the subgroup 
$P=\{g\in\G_M\mid e_kge_j=0\mbox{ if }1\le j<k\le n\}$
of $\G_M$,  
and define the mapping 
$$\psi\colon\U_M/\U_{M_0}\to\G_M/P,\quad 
u\U_{M_0}\mapsto uP.$$
Then $\psi$ is a real analytic $\U_M$-equivariant
diffeomorphism. 
\end{proposition}

\noindent\textbf{Proof.\  \  }
The mapping $\psi$ is clearly real analytic. 
Next note that $\U_M\cap P=\U_{M_0}$, hence $\psi$ is
injective. 

The fact that $\psi$ is surjective can be equivalently
expressed 
by the following assertion: 
For all $g\in\G_M$ there exist $u\in\U_M$ and $q\in P$
such that 
$g=uq$. 
In order to prove this fact, denote $p_0=0$ and 
$p_j=e_1+\cdots+e_j$ for $j=1,\dots,n$, 
and note that 
$$P=\{g\in\G_M\mid gp_j=p_jgp_j\mbox{ for
}j=1,\dots,p_n\}.$$
We have 
$$p_1\le p_2\le\cdots\le p_n
\mbox{ and
}\lf(gp_1)\le\lf(gp_2)\le\cdots\le\lf(gp_n).$$ 
Now let us fix $j\in\{1,\dots,n\}$.
Recall that for each $b \in M $ one denotes by $\lf(b)$
the smallest projection $p \in M $  satisfying $pb = b$.
By Lemma~3.2(iii) in
\cite{Be02} we have 
$\lf(gp_j)\sim p_j$ and $\lf(gp_{j-1})\sim p_{j-1}$. 
Since $p_j$ is a finite projection,
it then follows that 
$\lf(gp_j)-\lf(gp_{j-1})\sim p_j-p_{j-1}=e_j$ 
(see for instance Exercise~6.9.8 in \cite{KR97}). 
In other words, there exists $v_j\in M$ such that 
$v_jv_j^*=e_j$ and $v_j^*v_j=\lf(gp_j)-\lf(gp_{j-1})$.

Now, denoting 
$$u=v_1+\cdots+v_n\in M$$
it is easy to check that $u\in\U_M$ 
and $u^*e_ju=\lf(gp_j)-\lf(gp_{j-1})$ for
$j=1,\dots,n$. 
By summing up these equalities for $j=1,\dots,k$, 
we get 
$u^{-1}\lf(gp_k)=p_ku^*$ for each $k\in\{1,\dots,n\}$.

It then follows that for every $k$ we have 
$\lf(u^{-1}\lf(gp_k))=\lf(p_ku^*)$, 
whence 
$\lf(u^{-1}gp_k)=p_k$ by Lemma~3.2(i) in~\cite{Be02}. 
Consequently $q:=u^{-1}g\in P$, and we are done. 
\quad $\blacksquare$

\section{Homogeneous vector bundles}

The present section is devoted to a general
construction of homogeneous vector bundles in a
$C^*$-algebraic setting. 

\begin{construction}\label{homog}
\normalfont
Let $A$ be a unital $C^*$-algebra,  
$B$ a unital sub-$C^*$-algebra of $A$ and 
$\varphi\colon A\to{\mathbb C}$ a state such that
there exists a 
conditional expectation $E\colon A\to B$ with
$\varphi\circ E=\varphi$. 
By {\bfi conditional expectation\/} we mean that $E$ is a
bounded linear mapping such that 
$\|E\|=1$ and $E$ is idempotent, that is, $E^2=E$. 
It then follows by the theorem of Tomiyama 
(see \cite{To57}~or~\cite{Sak71}) 
that $E$ has the additional properties: 
\begin{gather}
E(a^*)=E(a)^*,\label{self-adj} \\
0\le E(a)^*E(a)\le E(a^*a),\label{positive} \\
E(b_1ab_2)=b_1E(a)b_2,\label{bimodule}
\end{gather}
for all $a\in A$ and $b_1,b_2\in B$. 

Now let 
$$\rho\colon A\to\Bc(\Hc)\quad\mbox{ and }\quad
\rho_\varphi\colon B\to\Bc(\Hc_\varphi)$$ 
be the GNS 
unital $*$-representations of $A$ and $B$
corresponding to 
$\varphi$ and  $\varphi|_B$, respectively.
We recall that, for instance, $\Hc$ is the Hilbert
space obtained from 
$A$ by factoring out the null-space of 
the nonnegative definite Hermitian sesquilinear form 
$$A\times A\ni(a_1,a_2)\mapsto\varphi(a_2^*a_1) \in
\mathbb{C}$$ 
and then taking the completion, and for
each
$a\in A$ the operator 
$\rho(a)\colon\Hc\to\Hc$ is the one obtained from 
the linear mapping $a'\mapsto aa'$ on $A$. 
The Hilbert space $\Hc_\varphi$ is obtained similarly,
using the restriction of the aforementioned
sesquilinear form to $B$. 
It is easy to see that 
$\Hc_\varphi$ is a closed subspace of $\Hc$ and 
for all $b\in B$ 
we have a commutative diagram 
\begin{equation}\label{GNS}
\begin{CD}
A \to \Hc @>{\rho(b)}>> \Hc\\
\hskip-25pt@V{E}VV \hskip-110pt@VV{P_{\Hc_\varphi}}V
\hskip-10pt 
@VV{P_{\Hc_\varphi}}V\\ 
B \to \Hc_\varphi @>{\rho_\varphi(b)}>> \Hc_\varphi
\end{CD}
\end{equation}
where $P_{\Hc_\varphi}\colon\Hc\to\Hc_\varphi$ 
denotes the orthogonal projection, 
while $A\to\Hc$ and $B \to \Hc_\varphi$ are the
inclusions. 
The fact that $E\colon A\to B$ 
is continuous with respect to the GNS scalar products
on $A$ and $B$ 
(whence its ``extension'' by continuity to
$P_{\Hc_\varphi}$) 
follows since for all $a\in A$ we have 
by inequality~\eqref{positive} that 
$E(a)^*E(a)\le E(a^*a)$, 
hence $\varphi(E(a)^*E(a))
\le\varphi(E(a^*a))=\varphi(a^*a)$.

By restriction we get a norm-continuous 
unitary representation 
of the unitary group of $B$, 
$$\rho|_{\U_{B}}\colon\U_{B}\to\U_{\Bc(\Hc_\varphi)}.$$
Since $\U_B$ is a Lie group with the topology
inherited from $\U_A$
and the self-adjoint mapping $E$ gives a continuous
projection 
of the Lie algebra of $\U_A$ onto the Lie algebra of
$\U_B$, 
it then follows by Proposition~11 in Chapter~III, \S
1.6 in 
\cite{Bo72} that
we have a principal $\U_B$-bundle 
$$\pi_B\colon\U_A\to\U_A/\U_B, \quad
u\mapsto u\cdot\U_B.$$
Now we can construct  
the $\U_A$-homogeneous vector bundle associated with 
$\pi_B$ and $\rho|_{\U_B}$, which we denote by 
$$\Pi_{\varphi,B}\colon D_{\varphi,B}\to\U_A/\U_B$$
(see for instance 6.5.1 in~\cite{Bo67} for this
construction). 
We recall that 
$D_{\varphi,B}=\U_A\times_{\U_B}\Hc_{\varphi}$, 
in the sense that
$D_{\varphi,B}$ is the quotient of 
$\U_A\times\Hc_\varphi$ by the equivalence relation
$\sim$ defined by
$$(u_1,f_1)\sim(u_2,f_2)\iff (\exists v\in\U_B)\quad 
u_1=u_2v\mbox{ and }f_1=\rho_\varphi(v^{-1})f_2,$$ 
while the mapping $\Pi_{\varphi,B}$ takes the
equivalence class 
$[(u,f)]$ 
of any pair $(u,f)$ to $u\cdot\U_B$. 
\end{construction}

\begin{definition}\label{Definition6}
\normalfont
With the notations of Construction~\ref{homog}, 
the $\U_A$-homogeneous vector
bundle $\Pi_{\varphi,B}\colon
D_{\varphi,B} \rightarrow \U_A/\U_B$
is called the {\bfi homogeneous  bundle
associated with $\varphi$ and $B$}. 
\end{definition}

\begin{remark}\label{Remark7}
\normalfont
In the setting of Construction~\ref{homog}, there
exist additional 
structures on the vector bundle 
$\Pi_{\varphi,B}\colon D_{\varphi,B}\to\U_A/\U_B$,
which we discuss below.
\begin{itemize}
\item[{\rm(i)}] Denote by 
$\alpha_B\colon\U_A\times\U_A/\U_B
\to\U_A/\U_B$ 
the natural action of $\U_A$ on $\U_A/\U_B$. 
Define 
the mapping  
$$\beta_{\varphi,B}\colon\U_A\times D_{\varphi,B}\to
D_{\varphi,B}$$
that takes the pair consisting of $u'\in\U_A$ and the
equivalence class of 
$(u,f)\in\U_A\times\Hc_\varphi$ to the equivalence
class of 
$(u'u,f)$. 
Then $\beta_{\varphi,B}$ is a real analytic action of 
$\U_A$ on $D_{\varphi,B}$ and the diagram 
$$\begin{CD}
\U_A\times D_{\varphi,B} @>{\beta_{\varphi,B}}>>
D_{\varphi,B} \cr
@V{\id_{\U_A}\times\Pi_{\varphi,B}}VV
@VV{\Pi_{\varphi,B}}V \cr
\U_A\times\U_A/\U_B @>{\alpha_B}>> \U_A/\U_B
\end{CD}$$
is commutative.

\item[{\rm(ii)}] The bundle $\Pi_{\varphi,B}$ is
actually a Hermitian 
vector bundle, in the sense that its fibers come
equipped 
with structures of 
complex Hilbert spaces, and moreover for each
$u\in\U_A$ 
the mapping 
$\beta_{\varphi,B}(u,\cdot)\colon D_{\varphi,B}\to
D_{\varphi,B}$ 
is (bounded linear and) fiberwise unitary.  
\end{itemize}
\end{remark}

\begin{example}\label{extremes}
\normalfont
Let $A$ be a unital $C^*$-algebra and 
$\varphi\colon A\to{\mathbb C}$ a state with the
corresponding GNS representation 
$\rho\colon A\to\Bc(\Hc)$. 
We describe below the extreme situations of the unital
sub-$C^*$-algebra $B$ of $A$, when the hypothesis of 
Construction~\ref{homog} 
is satisfied. 
\begin{itemize}
\item[{\rm(i)}] For $B=A$ we can take $E=\id_A$. 
In this case the homogeneous vector bundle associated
with 
$\varphi$ and $B$ is just the vector bundle with the
base 
reduced to a single point and with the fiber equal to
$\Hc$. 
\item[{\rm(ii)}] The other extreme situation is for 
$B={\mathbb C}\1$, when we can take
$E(\cdot)=\varphi(\cdot)\1$. 
Then, denoting ${\mathbb T}=\{z\in{\mathbb
C}\mid|z|=1\}$, 
it follows that 
the homogeneous vector bundle associated with 
$\varphi$ and $B$ is a line bundle whose base is 
the Banach-Lie group 
$\U_A/{\mathbb T}\1$.  
\end{itemize}
\end{example}

\begin{example}\label{W*}
\normalfont
If $M$ is a $W^*$-algebra and $0\le\varphi\in
M_*$, then there always exists a conditional
expectation
$$E_\varphi\colon M\to M^\varphi\mbox{ with
}\varphi\circ E_\varphi=\varphi.$$
We recall that the existence of $E_\varphi$ follows by
Theorem~4.2 in Chapter~IX of~\cite{Ta03} 
in the case when $\varphi$ is faithful. 
In the general case, 
denote by $p=p^*=p^2\in M$ the support of $\varphi$. 
Since in a $W^*$-algebra every element is a linear
combination 
of unitary elements, 
it follows by Lemma~2.7 in~\cite{BR04} 
that, denoting
$\varphi_p:=\varphi|_{(pMp)}\in(pMp)_*$, we have 
$$M^\varphi=\{a\in M\mid ap=pa\mbox{ and }pap\in(pMp)^{\varphi_p}\}
=(pMp)^{\varphi_p}\oplus(\1-p)M(\1-p).$$
Now, since $\varphi_p$ is faithful on $pMp$, 
it follows by the aforementioned theorem from~\cite{Ta03}
that there exists a conditional expectation 
$E_{\varphi_p}\colon pMp\to(pMp)^{\varphi_p}$ 
with $\varphi_p\circ E_{\varphi_p}=E_{\varphi_p}$. 
Then 
$$E_\varphi\colon M\to M^{\varphi},\quad 
a\mapsto E_{\varphi_p}(pap)+(\1-p)a(\1-p)$$
is a conditional expectation with $\varphi\circ
E_\varphi=\varphi$ 
(see also Remark~2.5 in \cite{BR04}). 
\end{example}

\begin{example}\label{C*}
\normalfont
Let $A$ be a unital $C^*$-algebra,  
$B$ a unital sub-$C^*$-algebra of $A$ and 
$E\colon A\to B$ a conditional expectation 
as in Construction~\ref{homog}. 
If $\varphi_0\colon B\to{\mathbb C}$ is a
\textit{pure} 
state of $B$, then it is easy to see that 
$\varphi:=\varphi_0\circ E\colon A\to{\mathbb C}$ 
is in turn a pure state of $A$
provided it is the unique extension of $\varphi_0$ to $A$.
Since pure states lead to irreducible representations,
this easy remark can be viewed as a version of
Theorem~2.5 
in~\cite{BH98} asserting that (roughly speaking)
irreducibility of the isotropy representation implies
global irreducibility. 
(Compare also Theorem~\ref{Theorem12} below.) 

We refer to the papers \cite{Ar99}~and~\cite{Ar01} 
for situations when, converse to the above remark, 
the existence of a unique conditional expectation from
$A$ onto $B$ 
follows from the assumption that each pure state of
$B$ extends uniquely to a pure state of $A$. 
\end{example}

\section{Reproducing kernels}

In what follows, 
by {\bfi continuous vector bundle} we mean a continuous
mapping 
$\Pi\colon D\to T$, where  
$D$ and $T$ are topological spaces, 
the fibers $D_t:=\Pi^{-1}(t)$ ($t\in T$) are complex
Banach spaces, 
and 
$\Pi$ is locally trivial with fiberwise bounded linear
trivializations. 
We say that $\Pi$ is {\bfi Hermitian} if it is equipped
with 
a continuous Hermitian structure 
$$(\cdot\mid\cdot)\colon D\oplus D\to{\mathbb C}$$
that makes each fiber $D_t$ into a complex Hilbert
space. 
We define in a similar manner the {\bfi smooth} vector
bundles and 
{\bfi holomorphic} vector bundles. 
For instance, for a Hermitian holomorphic vector
bundle we require that 
both $D$ and $T$ are complex Banach manifolds, $\Pi$
is holomorphic and 
the Hermitian structure is smooth. 
We denote by $\Cc(T,E)$, $\Cc^\infty(T,E)$ and
$\Oc(T,E)$ 
the spaces of continuous, $\Cc^\infty$ and holomorphic
sections of the bundle $\Pi\colon T\to E$, 
respectively, whenever the corresponding smoothness
condition makes sense. 
We refer to \cite{La01}~and~\cite{AMR88} for basic
facts on vector bundles.

\begin{definition}\label{Definition8}
\normalfont
Let $\Pi\colon D\to T$ be a continuous Hermitian
vector bundle,
$p_1,p_2\colon T\times T\to T$ the natural projections, 
and 
$$p_j^*\Pi\colon p_j^*D\to T\times T$$
the pull-back of $\Pi$ by $p_j$, $j=1,2$,
which are in turn
Hermitian vector bundles. 
Also consider the continuous vector bundle 
$\Hom(p_2^*\Pi,p_1^*\Pi)$ over $T\times T$, 
whose fiber over $(s,t)\in T\times T$ 
is the Banach space $\Bc(D_t,D_s)$. 

A {\bfi positive definite reproducing kernel} on
$\Pi$ 
is a continuous section 
$$K\in\Cc(T\times T,\Hom(p_2^*\Pi,p_1^*\Pi))$$
having the property that
for every integer $n\ge1$ and all choices of 
$t_1,\dots,t_n\in T$ and $\xi_1\in
D_{t_1},\dots,\xi_n\in D_{t_n}$, we have 
$$\sum_{j,l=1}^n(K(t_l,t_j)\xi_j\mid\xi_l)_{D_{t_l}}\ge0.$$
If $\Pi$ is a holomorphic vector bundle,
we say that the reproducing kernel $K$ is {\bfi
holomorphic} 
if for each $t\in T$ and every $\xi\in D_t$ 
we have $K(\cdot,t)\xi\in\Oc(T,D)$.  \
\medskip

The following properties of a reproducing kernel $K$ 
are immediate consequences of the above definition:
\begin{itemize}
\item[{\rm(i)}] For all $t\in T$ 
we have $K(t,t)\ge0$ in $\Bc(D_t)$. 

\item[{\rm(ii)}] For all $t,s\in T$ we have 
$K(t,s)\in\Bc(D_s,D_t)$, $K(s,t)\in\Bc(D_t,D_s)$ and 
$K(t,s)^*=K(s,t)$. 
\end{itemize}
\end{definition}

The first part of the next statement is a version of 
Theorem~1.4 in \cite{BH98}. 
The proof consists in adapting to the present setting 
a number of basic ideas from the classical 
theory of reproducing kernel Hilbert spaces 
(see for instance~\cite{Ha82} and \cite{CG04}).

\begin{theorem}\label{Theorem9}
Let $\Pi\colon D\to T$ be a continuous Hermitian
vector bundle, 
denote by 
$p_1,p_2\colon T\times T\to T$ the projections and 
consider
$K\in\Cc(T\times T,\Hom(p_2^*\Pi,p_1^*\Pi))$. 
Then $K$ is a reproducing kernel on $\Pi$ if and only
if 
there exists a linear mapping
$\iota\colon\Hc\to\Cc(T,D)$,
where $\Hc$ is a complex Hilbert space 
and $\ev^\iota_t:=(\iota(\cdot))(t)\colon\Hc\to D_t$ 
is a bounded linear operator for all $t\in T$, such
that 
\begin{equation}
\label{k in terms of ev}
K(s,t)=\ev^\iota_s
(\ev^\iota_t)^*\colon D_t\to D_s,
\end{equation}
for all $t\in T$. 
If this is the case, then the mapping $\iota$ can be
chosen to be injective. 
If, moreover, $\Pi$ is a holomorphic vector bundle then 
the reproducing kernel $K$ is holomorphic if and only
if 
any mapping $\iota$ as above has the range contained
in 
$\Oc(T,D)$. 
\end{theorem}

\noindent\textbf{Proof.\  \  }
Given a mapping $\iota$ as described in the statement,
it is clear that $K$ defined by \eqref{k in terms of
ev} is a reproducing kernel. 

Conversely, given the reproducing kernel $K$, 
consider the complex linear subspace $\Hc^K_0 \subset
\Cc(T,D)$  spanned by
$\{K_\xi:=K(\cdot,\Pi(\xi))\xi\mid\xi\in D\}$ 
 and denote by $\Hc^K$ its completion with respect to
the scalar product 
defined by 
$$(\Theta\mid\Delta)_{\Hc^K}
=\sum_{j,l=1}^n(K(\Pi(\xi_l),\Pi(\eta_j))\eta_j\mid 
\xi_l)_{D_{\Pi(\xi_l)}}, $$
where  $\Theta=\sum\limits_{j=1}^n
K(\cdot,\Pi(\eta_j))\eta_j,\,
\Delta=\sum\limits_{l=1}^n K(\cdot,\Pi(\xi_l))\xi_l 
\in \Hc^K_0$.
In order to see that this formula is independent
on the choices used to define $\Theta $ and
$\Delta$,  let us first
note the reproducing kernel property
\begin{equation}\label{reprod}
(\Theta\mid K_{\xi})_{\Hc^K}
=(\Theta(\Pi(\xi))\mid\xi)_{D_{\Pi(\xi)}}
\end{equation}
which holds whenever $\xi\in D$ and $\Theta\in\Hc^K_0$.
Now, for $\Theta$ and $\Delta$ as above we get 
$$(\Theta\mid\Delta)_{\Hc^K}
=\sum_{l=1}^n(\Theta(\Pi(\xi_l))\mid\xi_l)_{D_{\Pi(\xi_l)}}.$$
Hence if $\Theta=0$, then
$(\Theta\mid\Delta)_{\Hc^K}=0$.  

Moreover, 
$$(\Theta\mid\Theta)_{\Hc^K}
=\sum_{j,l=1}^n(K(\Pi(\eta_l),\Pi(\eta_j))\eta_j\mid 
\eta_l)_{D_{\Pi(\eta_l)}}\ge0,$$
which shows that $(\cdot\mid\cdot)_{\Hc^K}$ is a nonnegative
definite Hermitian sesquilinear form on $\Hc^K_0$. 
Therefore, it satisfies the Cauchy-Schwarz inequality and,
in particular, the reproducing kernel property~\eqref{reprod}
implies that 
$$|(\Theta(\Pi(\xi))\mid\xi)_{D_{\Pi(\xi)}}|
\le(\Theta\mid\Theta)_{\Hc^K}^{1/2}
(K_{\xi}\mid K_{\xi})_{\Hc^K}^{1/2}
\quad\mbox{ for all }\xi\in D.$$
Thus, if $(\Theta\mid\Theta)_{\Hc^K}=0$, then
$\Theta=0$.  
On the other hand, for each $\xi\in D$ 
we  have 
$$
(K_{\xi}\mid K_{\xi})_{\Hc^K}
=(K(\Pi(\xi),\Pi(\xi))\xi\mid\xi)_{D_{\Pi(\xi)}}
\le\|K(\Pi(\xi),\Pi(\xi))\|\|\xi\|_{D_{\Pi(\xi)}}^2.
$$
Now it follows by this inequality along with the
previous one 
that for all $\Theta\in\Hc^K_0$ and $t\in T$ we have
$$\|\Theta(t)\|_{D_t}\le\|K(t,t)\|^{1/2}
(\Theta\mid\Theta)_{\Hc^K}^{1/2}
=\|K(t,t)\|^{1/2}
\|\Theta\|_{\Hc^K}.$$
Consequently we can uniquely define for each $t\in T$ 
the bounded linear mapping 
$\ev^\iota_t\colon\Hc^K\to D_t$ 
satisfying 
$\ev^{\iota}_t(\Theta)=\Theta(t)$ for all  
$\Theta\in\Hc^K_0$. 
Then we define $\iota\colon\Hc^K\to\Cc(T,D)$  by
$(\iota(\cdot))(t): = \ev^\iota_t(\cdot)\colon\Hc^K\to
D_t$  for all $t\in T$.  
Now note that the reproducing kernel
property~\eqref{reprod}  extends by continuity to arbitrary 
$\Theta\in\Hc^K$ under the following
form: 
\begin{equation}\label{REPROD}
(\Theta\mid K_{\xi})_{\Hc^K}
=((\iota(\Theta))(\Pi(\xi))\mid\xi)_{D_{\Pi(\xi)}}
\mbox{ for all }\Theta\in\Hc^K
\mbox{ and }\xi\in D.
\end{equation}
In particular, \eqref{REPROD} shows that
$\iota\colon\Hc^K\to\Cc(T,D)$ 
is injective: if $\iota(\Theta)=0$, then $\Theta\perp
K_\xi$ 
in $\Hc^K$ for all $\xi\in D$, 
hence $\Theta=0$ since the linear span of
$\{K_\xi\mid\xi\in D\}$ 
is dense in $\Hc^K$. 
The proof of the first part of the statement is
finished. 

Now assume that $\Pi$ is a holomorphic vector bundle
and let $\iota$ be a corresponding injective mapping. 
Note that the expression \eqref{k in terms of ev} of $K$
in terms of the mappings 
$\ev_s^\iota$ implies that for each $t\in T$ and
$\xi\in D_t$ we have 
\begin{equation}\label{formula}
K(s,t)\xi=\ev_s^\iota((\ev_t^\iota)^*\xi)
=\iota((\ev_t^\iota)^*\xi)(s)
\quad\mbox{ for all }s\in T.
\end{equation}
Thus, if $\Ran\iota\subseteq\Oc(T,D)$ then $K$ is a
holomorphic kernel. 
Conversely, assume that $\iota$ is injective and the
kernel $K$ is holomorphic. 
Then \eqref{formula} shows that
$\iota(h)\in\Oc(T,D)$ 
for all 
$h\in\Hc_0:=\bigcup\limits_{\xi\in
D}\Ran(\ev_{\Pi(\xi)}^\iota)^*$ 
($\subseteq\Hc$). 
Note that in $\Hc$ we have 
$$\Hc_0^\perp
=\bigcap\limits_{\xi\in
D}(\Ran(\ev_{\Pi(\xi)}^\iota)^*)^\perp
=\bigcap\limits_{\xi\in D}\Ker(\ev_{\Pi(\xi)}^\iota)
=\{0\},$$
where the latter equality follows since $\iota$ is
injective. 
Consequently we know that $\iota(h)\in\Oc(D,T)$ when
$h$ runs 
through a dense linear subspace of $\Hc$. 
Now note that, as above, we can show that 
$$\|(\iota(h))(t)\|_{D_t}\le\|K(t,t)\|\,\|h\|_{\Hc}
\quad\mbox{ for all }h\in\Hc,$$ 
hence $\iota\colon\Hc\to\Cc(D,T)$ is continuous when 
$\Cc(D,T)$ is equipped with the topology of uniform
convergence on the 
compact subsets of $D$. 
On the other hand, $\Oc(D,T)$ is a closed subspace of
$\Cc(D,T)$ 
with respect to that topology 
(see e.g., Corollary~1.14 in~\cite{Up85})). 
Since we have already seen that
$\iota(\Hc_0)\subseteq\Oc(D,T)$ 
and $\Hc_0$ is dense in $\Hc$, it then follows that 
$\Ran\iota\subseteq\Oc(D,T)$, as desired.
\quad $\blacksquare$

\section{Reproducing kernels and GNS representations}

In this section we obtain our main results concerning
geometric realizations of GNS representations in
spaces of sections of Hermitian vector bundles 
(see Theorems \ref{Theorem12}~and~\ref{holomorphy}). 
Since the corresponding spaces of sections will be
reproducing kernel Hilbert spaces, we begin by
constructing the reproducing kernels.
In the extreme case $B={\mathbb C}\1$
(see Example~\ref{extremes}(ii))
the kernel $K_{\varphi,B}$ from the following construction is related
to a certain hermitian kernel studied in \cite{CG01}
(more precisely, see formula~(6.2) in \cite{CG01}).

\begin{construction}\label{kernel}
\normalfont
Let $A$ be a unital $C^*$-algebra,  
$B$ a unital sub-$C^*$-algebra of $A$ and 
$\varphi\colon A\to{\mathbb C}$ a state such that
there exists a 
conditional expectation $E\colon A\to B$ with
$\varphi\circ E=\varphi$.
Recall from Definition~\ref{Definition6} the
homogeneous vector bundle 
$$\Pi_{\varphi,B}\colon
D_{\varphi,B}=\U_A\times_{\U_{B}}\Hc_{\varphi}
\to\U_A/\U_{B}$$
associated with $\varphi$ and $B$. 
Then we define 
$$\iota_{\varphi,B}\colon\Hc\to\Cc(\U_A/\U_B,D_{\varphi,B}),
\quad \text{by} \quad 
\iota_{\varphi,B}(h)(u\U_B)=[(u,P_{\Hc_\varphi}(\rho(u)^{-1}h))]$$
for all $h\in\Hc$ and $u\in\U_A$. 
We also define 
$$K_{\varphi,B}\in\Cc\bigl(\U_A/\U_{B}\times\U_A/\U_{B},
\Hom(p_2^*(\Pi_{\varphi,B}),p_1^*(\Pi_{\varphi,B}))\bigr)$$
by 
$$K_{\varphi,B}(u_1\U_{B},u_2\U_{B})[(u_2,f)]
=[(u_1,P_{\Hc_\varphi}(\rho(u_1^{-1}u_2)f))],$$
where $u_1,u_2\in\U_A$,  
$f\in\Hc_\varphi$, and 
$p_1,p_2\colon\U_A/\U_{B}\times\U_A/\U_{B}
\to\U_A/\U_{B}$ 
are the projections. 

We note that both $\iota_{\varphi,B}$ and
$K_{\varphi,B}$ 
are well defined by the commutativity of the
diagram~\eqref{GNS}. 
\end{construction}

\begin{definition}\label{Definition10}
\normalfont
In the setting of Construction~\ref{kernel}, 
the operator 
$\iota_{\varphi,B}\colon\Hc\to\Cc(\U_A/\U_{B},D_{\varphi,B})$
is called the {\bfi realization operator} 
associated with $\varphi$ and $B$ and the map
$K_{\varphi,B}\in\Cc\bigl(\U_A/\U_{B}\times\U_A/\U_{B},
\Hom(p_2^*(\Pi_{\varphi,B}),p_1^*(\Pi_{\varphi,B}))\bigr)$
is called 
the {\bfi reproducing kernel} associated with $\varphi$
and $B$. 
\end{definition}

\begin{remark}\label{Remark11}
\normalfont
\begin{itemize}
\item[{\rm(i)}]
It will follow from the proof of
Theorem~\ref{Theorem12} below 
that the reproducing kernel $K_{\varphi,B}$ associated
with a 
state $\varphi$ and $B$ is 
indeed a reproducing kernel in the sense of 
Definition~\ref{Definition8}.

\item[{\rm(ii)}]
Recall the action $\beta_{\varphi,B}\colon\U_A\times
D_{\varphi,B}
\to D_{\varphi,B}$ 
from Remark~\ref{Remark7}(i) and, for all
$u,u'\in\U_A$ and $f\in\Hc_\varphi$,
denote 
$u'\cdot
[(u,f)]:=\beta_{\varphi,B}(u',[(u,f)])=[(u'u,f)]$  for
the sake of simplicity.  Then, in the notation of 
Definition~\ref{Definition10}, for all
$u,u_1,u_2\in\U_A$ 
and $f\in\Hc$ we have 
$$\begin{aligned}
K_{\varphi,B}(u\cdot u_1\U_{B},u_2\U_{B})[(u_2,f)]
&=[(uu_1,P_{\Hc_\varphi}(\rho(u_1^*u^*u_2)f))] \\
&=u\cdot [(u_1,P_{\Hc_\varphi}(\rho(u_1^*u^*u_2)f))] \\
&=u\cdot
K_{\varphi,B}(u_1\U_{B},u^*u_2\U_{B})[(u^*u_2,f)]\\
&=u\cdot K_{\varphi,B}(u_1\U_{B},u^*u_2\U_{B})
u^*\cdot [(u_2,f)].
\end{aligned}$$
Thus
$$u^*\cdot K_{\varphi,B}(u\cdot
t_1,t_2)\xi=K_{\varphi,B}(t_1,u^*\cdot t_2)(u^*\cdot
\xi)$$ for $u\in\U_A$,
$t_1,t_2\in\U_A/\U_{B}$  and $\xi\in
(D_{\varphi,B})_{t_2}$.
\end{itemize}
\end{remark}

In the statement of the next theorem we refer to the
correspondence 
established in Theorem~\ref{Theorem9}. 
Representations of this type were also considered from
another point of view in~\cite{Ne02}. 

\begin{theorem}\label{Theorem12}
Let $A$ be a unital $C^*$-algebra, 
$B$ a unital sub-$C^*$-algebra of $A$ and 
$\varphi\colon A\to{\mathbb C}$ a state such that
there exists a 
conditional expectation $E\colon A\to B$ with
$\varphi\circ E=\varphi$.
Consider the GNS representation 
$\rho\colon A\to\Bc(\Hc)$ 
associated with $\varphi$. 
Then the realization operator 
$\iota_{\varphi,B}\colon\Hc\to\Cc(\U_A/\U_{B},D_{\varphi,B})$
associated with $\varphi$ and $B$ is injective, 
has the property that 
$(\iota_{\varphi,B}(\cdot))(s)\colon\Hc\to(D_{\varphi,B})_s$
is bounded linear for all $s\in\U_A/\U_{B}$, 
and the corresponding reproducing kernel is precisely 
the reproducing kernel $K_{\varphi,B}$ associated with
$\varphi$ and $B$. 
Moreover, $\iota_{\varphi,B}$ intertwines the unitary
representation 
$\rho|_{\U_A}$ of $\U_A$ on $\Hc$ 
and the natural representation of $\U_A$ 
by linear mappings on
$\Cc(\U_A/\U_{B},D_{\varphi,B})$. 
\end{theorem}

\noindent\textbf{Proof.\  \  }
According to Theorem~\ref{Theorem9}, 
the asserted relationship between $\iota_{\varphi,B}$
and $K_{\varphi,B}$ 
will follow as soon as we 
prove that for all $s,t\in\U_A/\U_{B}$ 
the equality 
$$K_{\varphi,B}(s,t)=\ev_s^{\iota_{\varphi,B}}(\ev_t^{\iota_{\varphi,B}})^*
\colon(D_{\varphi,B})_t\to(D_{\varphi,B})_s$$
holds, where 
$$\ev_s^{\iota_{\varphi,B}}
=(\iota_{\varphi,B}(\cdot))(s)\colon\Hc\to(D_{\varphi,B})_s$$
and similarly for $\ev_t^{\iota_{\varphi,B}}$. 

To this end, let $u_1,u_2\in\U_A$ such that 
$s=u_1\U_{B}$ and $t=u_2\U_{B}$. 
Then 
$$(D_{\varphi,B})_s=\left\{[(u_1,f)]\mid
f\in\Hc_\varphi\right\}$$
and 
$$\ev_s^{\iota_{\varphi,B}}(h):=[(u_1,P_{\Hc_\varphi}(\rho(u_1^{-1})h))]
=[(u_1,(P_{\Hc_\varphi}\circ\rho(u_1)^*)h))]
\quad\mbox{for all }h\in\Hc,$$
and a similar formula holds with $s$ replaced by $t$ 
and $u_1$ replaced by $u_2$. 
Now, 
since $(P_{\Hc_\varphi})^*$ is just the inclusion
mapping 
$\Hc_\varphi\hookrightarrow\Hc$, it follows that 
for an arbitrary element
$[(u_2,f)]\in(D_{\varphi,B})_t$ 
(where $f\in\Hc_\varphi$) we have 
$(\ev_t^{\iota_{\varphi,B}})^*[(u_2,f)]
=\rho(u_2)f$;
hence 
$$\ev_s^{\iota_{\varphi,B}}(\ev_t^{\iota_{\varphi,B}})^*[(u_2,f)]
=[(u_1,(P_{\Hc_\varphi}\circ\rho(u_1)^*)\rho(u_2)f))]
=K_{\varphi,B}(s,t)[(u_2,f)],$$
as desired. 

In order to prove that $\iota_{\varphi,B}$ is
injective, 
let $h\in\Hc$ 
such that $\iota_{\varphi,B}(h)=0$. 
By the definition of $\iota_{\varphi,B}$ 
(see Construction~\ref{kernel}), 
it follows that for all $u\in\U_A$ we have 
$P_{\Hc_\varphi}(\rho(u^{-1})h)=0$. 
On the other hand, $A$ is the linear span of $\U_A$. 
(For instance, if $a=a^*\in A$ and $\|a\|\le1$, 
then $a=(u+u^*)/2$, where
$u=a+\ie(\1-a^2)^{1/2}\in\U_A$.) 
Consequently $P_{\Hc_\varphi}(\rho(a)h)=0$ for all
$a\in M$. 
In particular, denoting by
$h_0\in\Hc_\varphi\subseteq\Hc$ 
the image of $\1\in A$ in $\Hc_\varphi$, we get 
$0=(\rho(a)h\mid h_0)=(h\mid\rho(a^*)h_0)$ for all
$a\in A$. 
Since $h_0$ is a cyclic vector for the representation 
$\rho\colon A\to\Bc(\Hc)$, according to the GNS
construction 
(see e.g., Proposition~1.100 and Definition~1.101 in
\cite{AS01}), 
it follows that $h\perp\Hc$, whence $h=0$. 
Thus $\iota_{\varphi,B}$ is injective. 

Since it is clear that 
$(\iota_\varphi(\cdot))(u \U_B)\colon\Hc\to
(D_{\varphi,B})_{u \U_B}$ is bounded linear for all
$u\in\U_A$, it follows  that the realization operator 
$\iota_{\varphi,B}\colon\Hc\to\Cc(\U_A/\U_{B},
D_{\varphi,B})$
is 
of the type occurring in Theorem~\ref{Theorem9}. 

It remains to prove the intertwining property of 
$\iota_{\varphi,B}$. 
So consider, as in the proof of 
Theorem~\ref{Theorem9},
the linear subspace $\Hc^{K_{\varphi,B}}_0$ 
spanned by 
$\{(K_{\varphi,B})_\xi:=K_{\varphi,B}(\cdot,\Pi_{\varphi,B}(\xi))\xi\mid
\xi\in D_{\varphi,B}\}$
in 
$\Cc(\U_A/\U_{B},D_{\varphi,B})$. 
Then it follows by 
Remark~\ref{Remark11}(ii) that for all $u\in\U_A$ we
have 
$u\cdot
(K_{\varphi,B})_\xi=(K_{\varphi,B})_{u^*\cdot\xi}$,
where the left-hand side denotes the action of $u$ on 
$(K_{\varphi,B})_\xi\in\Cc(\U_A/\U_{B},D_{\varphi,B})$.

In view of the construction of $\Hc^K$ in the proof of
Theorem~\ref{Theorem9}, 
it follows that the natural action of $u$ on 
$\Cc(\U_A/\U_{B},D_{\varphi,B})$ 
leaves $\Hc^{K_{\varphi,B}}$ invariant and actually
induces a unitary operator on it. 
In fact, 
for $t_1,\dots,t_n\in\U_A/\U_B$, 
$\xi_1\in(D_{\varphi,B})_{t_1},\dots,\xi_n\in(D_{\varphi,B})_{t_n}$
and 
$\Theta=\sum\limits_{j=1}^n(K_{\varphi,B})_{\xi_j}=\sum\limits_{j=1}^n
K_{\varphi,B}(\cdot,t_j)\xi_j$, 
we have 
$\|\Theta\|_{\Hc^{K_{\varphi,B}}}^2
=\sum\limits_{j,l=1}^n(K(t_l,t_j)\xi_j\mid\xi_l)_{(D_{\varphi,B})_{t_l}}$
and 
$$\begin{aligned}
\|u\cdot\Theta\|_{\Hc^{K_{\varphi,B}}}^2
&=\sum_{j,l=1}^n
(K(u^*\cdot t_l,h^*\cdot t_j)(u^*\cdot \xi_j)\mid
u^*\cdot \xi_l)_{(D_{\varphi,B})_{t_l}} \\
&=\sum_{j,l=1}^n
(u\cdot K(u^*\cdot t_l,h^*\cdot t_j)(u^*\cdot \xi_j)\mid
\xi_l)_{(D_{\varphi,B})_{t_l}} \\
&=\sum_{j,l=1}^n(K(t_l,t_j)\xi_j\mid\xi_l)_{(D_{\varphi,B})_{t_l}}
\qquad\qquad\qquad\qquad\qquad(\mbox{by
Remark~\ref{Remark11}(ii)})\\
&=\|\Theta\|_{\Hc^{K_{\varphi,B}}}^2.
\end{aligned}$$
The intertwining property of $\iota_{\varphi,B}$
is straightforward: 
for all $u,v\in\U_A$ and $h\in\Hc$ 
we have 
$$\begin{aligned}
\iota_{\varphi,B}(\rho(v)h)(u\U_{B})
&=[(u,P_{\Hc_\varphi}(\rho(u^{-1})\rho(v)h))] 
=[(u,P_{\Hc_\varphi}(\rho((v^{-1}u)^{-1})h))] \\
&=v\cdot
[(v^{-1}u,P_{\Hc_\varphi}(\rho((v^{-1}u)^{-1})h))]
=v\cdot \iota_{\varphi,B}(h)(v^{-1}u\U_{B})
\end{aligned}$$
as desired.
\quad $\blacksquare$

\begin{remark}\label{Remark13a}
\normalfont
As a by-product 
of the proof of Theorem~\ref{Theorem12} it follows  
that for each $f\in\Hc_\varphi$ we have 
$\iota_{\varphi,B}(f)=(K_{\varphi,B})(\cdot,\1\U_B)[(\1,f)]$, 
whence 
$\|\iota_{\varphi,B}(f)\|_{\Hc^{K_{\varphi,B}}}
=\|f\|_{\Hc_\varphi}$. 
\end{remark}

\begin{remark}\label{Remark13}
\normalfont
In the setting of Theorem~\ref{Theorem12}, if it
happens 
that $\varphi$ is a pure state then the fact that 
the realization operator 
$\iota_{\varphi,B}\colon\Hc\to\Cc(\U_A/\U_{B},D_{\varphi,B})$
is injective can be proved in an alternative way as
follows. Let $h\in H$ with $\iota_{\varphi,B}(h)=0$. 
Then for all $u\in\U_A$ we have 
$0=(\iota_{\varphi,B}(h))(u \U_B)
=[(u, P_{\Hc_\varphi}(\rho(u^{-1})h))]$. Since every
element of
$A$ is a linear combination of unitary elements,  we get
$P_{\Hc_\varphi}(\rho(A)h)=0$, 
hence the closure $\Hc_0$ of $\rho(A)h$ is contained
in 
$\Ker P_{\Hc_\varphi}$. 
On the other hand, we cannot have $P_{\Hc_\varphi}=0$,
since, denoting by $h_0$ the image of $\1\in B$ in
$\Hc_\varphi$, 
we have $P_{\Hc_\varphi}h_0=h_0\ne0$. 
Hence $\Hc_0\ne\Hc$. 
Now note that $\Hc_0$ is an invariant subspace for the
GNS representation 
$\rho\colon A\to\Bc(\Hc)$ and this representation is
irreducible since $\varphi$ is pure. 
Consequently we must have $\Hc_0=\{0\}$, whence $h=0$,
 and thus $\iota_{\varphi,B}$ is injective. 
\end{remark}

In connection with Remark~\ref{Remark13}, 
we note that pure normal states can exist only in the
case of 
$W^*$-algebras of type~I, for instance direct sums of 
algebras of the form $\Bc(\Hc)$ with $\Hc$ a Hilbert
space. 
We discuss below the case of finite-dimensional $\Hc$.

A few details on the infinite-dimensional case can be
found in 
\cite{CGM03}.

\begin{example}\label{Example14}
\normalfont
Let $n\ge1$ be an integer, $M=M_n({\mathbb C})$
with its structure of $W^*$-algebra, and think of 
${\mathbb C}^n$ as a Hilbert space with the usual
scalar product 
$(\cdot\mid\cdot)$. 
Next, denote 
$$h=\begin{pmatrix} 1\\ 0\\ \vdots \\ 0\end{pmatrix}
\mbox{ and }
p=\begin{pmatrix} 1 & 0 & \dots & 0 \\ 
           0 & 0 & \dots & 0 \\
           \vdots &\vdots &\ddots & \vdots \\
           0 & 0 & \dots & 0 \\ 
  \end{pmatrix}\in M.$$
Now define
$$\varphi\colon M\to{\mathbb C},\quad 
\varphi(a)=(ah\mid h)=\Tr(ap)=a_{11}
\mbox{ for }a=(a_{ij})_{1\le i,j\le n}\in M.$$ 
It is clear that $\varphi$ is a pure normal state of
$M$. 
We want to construct the GNS representation of $M$
with 
respect to $\varphi$ and to see what
Theorem~\ref{Theorem12} says in this special case. 

We have 
$$M_0:=\{a\in M\mid\varphi(a^*a)=0\}
=\{a\in M\mid ah=0\}
=\left\{\begin{pmatrix} 0 & * &\cdots & *\\
                 0 & * &\cdots & *\\
                 \vdots & \vdots & \ddots &\vdots \\
                 0 & * &\cdots & *\end{pmatrix}
\right\}$$
and 
$$\varphi(a^*a)=|a_{11}|^2+|a_{21}|^2+\cdots+|a_{n1}|^2
\quad\mbox{ for }
a=\begin{pmatrix} a_{11} & 0 &\cdots & 0\\
                   a_{21} & 0 &\cdots & 0\\
                   \vdots & \vdots & \ddots &\vdots \\
                   a_{n1} & 0 &\cdots & 0\end{pmatrix}.
$$
Hence the completion $\Hc$ of $M/M_0$ with respect to 
the scalar product induced by
$(a,b)\mapsto\varphi(b^*a)$
is just the Hilbert space ${\mathbb C}^n$ with the
usual scalar product, viewed as the set of column
vectors, 
and the natural mapping $M\to\Hc$ is 
$$\begin{pmatrix} a_{11} & a_{12} &\cdots & a_{1n}\\
                   a_{21} & a_{22} &\cdots & a_{2n}\\
                   \vdots & \vdots & \ddots &\vdots \\
                   a_{n1} & a_{n2} &\cdots &
a_{nn}\end{pmatrix}
\mapsto
\begin{pmatrix} a_{11}\\
         a_{21}\\
         \vdots\\
         a_{n1}\end{pmatrix}.$$ 
Moreover, according to the way matrices multiply, it
follows that for each $a\in M$ the operator 
$\rho(a)\colon\Hc\to\Hc$ that makes the diagram 
$$\begin{CD}
M @>{b\mapsto ab}>> M \\
@VVV @VVV \\
\Hc @>{\rho(a)}>> \Hc
\end{CD}$$
commutative is just the natural action of $a\in
M=M_n({\mathbb C})$ 
on ${\mathbb C}^n$. 
Thus the GNS representation of $M$ associated to
$\varphi$ is just 
the natural representation 
of $M=M_n({\mathbb C})$ on ${\mathbb C}^n$. 

Next, Proposition~\ref{Proposition2}(b) shows that 
$$%\begin{aligned}
M^\varphi
=\{p\}'=\left\{\begin{pmatrix} * & 0 &\cdots & 0\\
                 0 & * &\cdots & *\\
                 \vdots & \vdots & \ddots &\vdots \\
                 0 & * &\cdots & *\end{pmatrix}
\right\}
=\left\{\begin{pmatrix} z & 0\\ 0 & W\end{pmatrix}\,
\Big{|}\,
 z \in{\mathbb C},W\in M_{n-1}({\mathbb C})\right\} %\\
\simeq M_1({\mathbb C})\times M_{n-1}({\mathbb C}),
%\end{aligned}
$$
hence $M^\varphi\cap M_0=(\1-p)M(\1-p)\simeq
M_{n-1}({\mathbb C})$
and thus 
the space of the GNS representation 
$\rho_\varphi\colon M^\varphi\to\Bc(\Hc_\varphi)$ 
of $M^\varphi$ corresponding to the state 
$\varphi|_{M^\varphi}\colon M^\varphi\to{\mathbb C}$ 
is one-dimensional, i.e., $\Hc_\varphi={\mathbb C}$. 
Moreover $\rho_\varphi\begin{pmatrix} z & 0\\ 0 &
W\end{pmatrix}$ 
is the multiplication-by-$z$ operator on $\Hc_\varphi$
for all $z\in{\mathbb C}$ and $W\in M_{n-1}({\mathbb
C})$. 
Now  
$$\U_M/\U_{M^\varphi}=\U(n)/(\U(1)\times\U(n-1))
={\mathbb P}^{n-1}({\mathbb C}),$$
and it follows at once that the homogeneous vector
bundle 
$\Pi_{\varphi,M^\varphi}\colon 
D_{\varphi,M^\varphi}\to\U_M/\U_{M^\varphi}$ 
associated with $\varphi$ is dual to the tautological
line bundle over the complex projective space
${\mathbb P}^{n-1}({\mathbb C})$. 
Thus, in this special case, Theorem~\ref{Theorem12} says
that 
the natural representation of $\U(n)$ on ${\mathbb C}^n$ 
can be geometrically realized as a representation in
the finite-dimensional vector space of global
holomorphic sections of the dual to the tautological line bundle
over the $\U(n)$-homogeneous compact 
K\"ahler manifold ${\mathbb P}^{n-1}({\mathbb C})$, 
which is a special case of the Borel-Weil theorem
(see e.g., \cite{ES02} and \cite{Do97}). 
\end{example}

We now describe a situation when the homogeneous
vector bundle 
occurring in Theorem~\ref{Theorem12} is holomorphic
and the 
reproducing kernel Hilbert space consists only of
holomorphic sections. 

\begin{theorem}\label{holomorphy}
Let $M$ be a $W^*$-algebra with 
a faithful normal tracial state $\tau\colon
M\to{\mathbb C}$. 
Pick $a\in M$ such that $0\le a$, $\tau(a)=1$ and the
spectrum of $a$ is a finite set, and define 
$$\varphi\colon M\to{\mathbb C},\quad
\varphi(b)=\tau(ab).$$
Then the homogeneous vector bundle 
$\Pi_{\varphi,M^\varphi}\colon 
D_{\varphi,M^\varphi}\to\U_M/\U_{M^\varphi}$ 
associated with $\varphi$ and $M^\varphi$ is
holomorphic and the reproducing kernel
$K_{\varphi,M^\varphi}$ is holomorphic as well. 
Also, if 
$$\rho\colon M\to\Bc(\Hc)$$ 
is the GNS representation 
corresponding to the normal state $\varphi$, then 
the image of the realization operator 
$\iota_{\varphi,M^\varphi}\colon
\Hc\to\Cc(\U_M/\U_{M^\varphi},D_{\varphi,M^\varphi})$ 
consists of holomorphic sections.
\end{theorem}

\noindent\textbf{Proof.\  \  }
As in the proof of Proposition~\ref{Proposition4}, 
write $a=\lambda_1e_1+\cdots+\lambda_ne_n$, 
where $e_1,\dots,e_n\in M$ are orthogonal projections
such that 
$e_ie_j=0$ whenever $i\ne j$ and $e_1+\cdots+e_n=\1$, 
and observe that, by Proposition~\ref{Proposition2}(b),
we have 
$$M^\varphi=\{a\}'
=\{b\in M\mid be_j=e_jb\mbox{ for }j=1,\dots,n\}.$$
Next denote by 
$$\rho_\varphi\colon M_0\to\Bc(\Hc_\varphi)$$
the GNS representation corresponding to
$\varphi|_{M^\varphi}$. 
Then the representation
$\rho_\varphi|_{\G_{M^\varphi}}$ of the 
group $\G_{M^\varphi}$ can be extended to a
representation of 
the group 
$$P:=\{g\in\G_M\mid e_kge_j=0\mbox{ if }1\le j<k\le
n\}$$
by defining 
$$\widetilde{\rho}_\varphi\colon
P\to\Bc(\Hc_\varphi),\quad 
\widetilde{\rho}_\varphi(g):=\rho_\varphi(e_1ge_1+\cdots+e_nge_n).$$
Using this representation, it makes sense to consider
the space 
$\widetilde{D}_\varphi=\G_M\times_P\Hc_\varphi$,
that is, the quotient of $\G_M\times\Hc_\varphi$ 
by the equivalence relation given by 
$(gc,f)\sim(g,\widetilde{\rho}_\varphi(c)f)$ whenever
$c\in P$. 

For an arbitrary pair $(g,f)\in\G_M\times\Hc_\varphi$ 
we denote by $[(g,f)]\widetilde{\phantom{b}}$ its
equivalence class in 
$\widetilde{D}_\varphi$, in order to distinguish it
from the equivalence class 
$[(g,f)]\in D_{\varphi,M^\varphi}
=\U_M\times_{\U_{M^\varphi}}\Hc_\varphi$ when it
happens that 
$g\in\U_M$. 
Also we define 
$\widetilde{\Pi}_\varphi([(g,f)]\widetilde{\phantom{b}})
=gP\in\G_M/P$ for all 
$[(g,f)]\widetilde{\phantom{b}}\in\widetilde{D}_\varphi$,
and $\Psi\colon D_\varphi\to \widetilde{D}_\varphi$ 
by $[(u,f)]\mapsto[(u,f)]\widetilde{\phantom{b}}$. 
Then we get a commutative diagram 
$$\begin{CD}
D_\varphi @>{\Psi}>> \widetilde{D}_\varphi \\
@V{\Pi_{\varphi,M^\varphi}}VV
@VV{\widetilde{\Pi}_\varphi}V \\
\U_M/\U_{M^\varphi} @>{\psi}>> \G_M/P
\end{CD}$$
where $\psi\colon\U_M/\U_{M^\varphi}\to\G_M/P$ is the
$\U_M$-equivariant 
diffeomorphism given by
Proposition~\ref{Proposition16}. 
In this diagram, $\widetilde{\Pi}_\varphi$ is a
holomorphic vector bundle (according to its
construction) and $\Psi$ is an $\U_M$-equivariant real
analytic mapping that is fiberwise bounded linear. 
We will prove shortly  that $\Psi$ is bijective and
then, 
since $\psi$ is a real analytic diffeomorphism by 
Proposition~\ref{Proposition16}, it will follow that
the pair 
$(\Psi,\psi)$ gives an isomorphism of vector bundles. 
Thus $\Pi_{\varphi,M^\varphi}$ is a holomorphic vector
bundle and, moreover, 
for each $h\in\Hc$ we have 
$(\Psi\circ\iota_{\varphi,M^\varphi}(h)\circ\psi^{-1})(gP)
=[(g,P_{\Hc_\varphi}(\rho(g)^{-1}h))]\widetilde{\phantom{b}}$
for all $g\in G$, 
hence clearly 
$\Psi\circ\iota_{\varphi,M^\varphi}(h)\circ\psi^{-1}
\in\Oc(\G_M/P,\widetilde{D}_\varphi)$, 
whence 
$\iota_{\varphi,M^\varphi}\in\Oc(\U_M/\U_{M^\varphi},D_\varphi)$.

Thus it only remains to show that the mapping $\Psi$
is bijective. 
To see that it is injective, let $u_1,u_2\in\U_M$ and 
$f_1,f_2\in\Hc_\varphi$ such that 
$[(u_1,f_1)]\widetilde{\phantom{b}}
=[(u_2,f_2)]\widetilde{\phantom{b}}$. 
Then there exists $c\in P$ such that $u_1=u_2c$ and 
$f_1=\widetilde{\rho}_\varphi(c^{-1})f_2$. 
In particular, $c=u_2^{-1}u_1\in\U_M$, 
whence $c\in\U_M\cap P=\U_{M^\varphi}$, 
and then $[(u_1,f_1)]=[(u_2,f_2)]$. 
Thus $\Psi$ is injective. 
To prove that it is surjective, let 
$g\in\G_M$ and $f\in\Hc_\varphi$ arbitrary. 
According to the proof of
Proposition~\ref{Proposition16}, 
there exist $u\in\U_M$ and $c\in P$ such that $g=uc$, 
whence 
$[(g,f)]\widetilde{\phantom{b}}
=[(uc,f)]\widetilde{\phantom{b}}
=[(u,\widetilde{\rho}_\varphi(c)f)]\widetilde{\phantom{b}}
=\Psi([(u,\widetilde{\rho}_\varphi(c)f)])$. 
\quad $\blacksquare$

\medskip

\noindent {\bf Acknowledgments.} 
We thank M. Rieffel and H. Upmeier for several
discussions that motivated us to seek geometric
realizations of representations of Banach Lie  groups.
We also thank A.~Gheondea and B.~Prunaru for some useful comments
and for drawing our attention to certain pertinent references. 
The first author was partially supported by the Swiss NSF
through the SCOPES Program during a one month visit at
the EPFL; the excellent working conditions provided by
the EPFL are gratefully acknowledged. The second author
acknowledges the partial support of the  Swiss NSF.

\end{document}